\documentclass[12pt]{article}%
\usepackage[applemac]{inputenc}
\usepackage{amsmath,amssymb}
\usepackage{amsfonts}
\usepackage{amsmath,amssymb}
\usepackage[applemac]{inputenc}
\usepackage{amsmath,amssymb}
\usepackage{color}
\usepackage{color, colortbl}
\usepackage{amsmath}
\usepackage{amssymb}
\usepackage{graphicx}
\usepackage{tikz}
\usepackage{geometry}
\usetikzlibrary{arrows}
\usepackage{amsfonts}
\usepackage{amsmath,amssymb}
\usepackage[applemac]{inputenc}
\usepackage{amsmath,amssymb}
\usepackage{color}
\usepackage{amsmath}
\usepackage{amssymb}
\usepackage{graphicx}
\usepackage{tikz}
\newtheorem{theorem}{Theorem}[section]
\newtheorem{lemma}[theorem]{Lemma}

\newtheorem{definition}[theorem]{Definition}

\newtheorem{example}[theorem]{Example}

\newtheorem{remark}[theorem]{Remark}

\usetikzlibrary{arrows,shapes,automata,petri}
\newenvironment{myenumerate}{

\begin{enumerate}}{\end{enumerate}}
\newcommand{\dproof}{\noindent {Proof.} \quad}
\newcommand{\fproof}{\hfill $\square$ \bigskip}

\numberwithin{equation}{section}

\definecolor{LightCyan}{rgb}{0.88,1,1}

\def\RR{{\mathbb{ R}}}

\def\EE{{\mathbb{ E}}}

\def\1B{\text{1\!\!I}}

\def\<{\langle}
\def\>{\rangle}

\begin{document}

\title{The Donsker delta function and local time for McKean-Vlasov processes and applications}
\author{Nacira Agram$^{1}$ \& Bernt \O ksendal$^{2,*}$}
\date{23 February 2023}
\maketitle

\footnotetext[1]{Department of Mathematics, KTH Royal Institute of Technology 100 44, Stockholm, Sweden. \newline
Email: nacira@kth.se. Work supported by the Swedish Research Council grant (2020-04697).}

\footnotetext[2,*]{%
Department of Mathematics, University of Oslo, Norway. 
Email: oksendal@math.uio.no\\ (corresponding author).}

\begin{abstract}
The purpose of this paper is to establish a stochastic differential equation for the Donsker delta measure of the solution of a McKean-Vlasov (mean-field) stochastic differential equation.

If the Donsker delta measure is absolutely continuous with respect to Lebesgue measure, then its
Radon-Nikodym derivative is called the Donsker delta function. In that case it can be proved that the local time of such a process is simply the integral with respect to time of the Donsker delta function. Therefore we also get an equation for the local time of such a process.

For some particular McKean-Vlasov processes, we  find  explicit expressions for their Donsker delta functions and hence for their local times.
\end{abstract}

\textbf{Keywords :} Donsker delta function, local time, McKean-Vlasov process, Fokker-Planck equation

\textbf{MSC 2020 :} 60H15; 60H40; 60J35

\section{Introduction}

The Donsker delta function of a random variable or a stochastic process arises in many studies, including quantum mechanical particles on a circle \cite{LLSW}, financial markets with insider trading as in \cite{OR}, and in \cite{AO1} for financial markets with singular drift. It has also been used as a tool to determine explicit formulae for replicating portfolios in complete and incomplete markets, see \cite{DOP}.

Moreover, the Donsker delta function is also of interest because it can be regarded as a time-derivative of the local time. Therefore, explicit expressions for the Donsker delta function lead to explicit formulae of the local time. 

For example, if we let $B$ be a Brownian motion defined on a filtered probability space  $(\Omega,\mathcal{F},\mathbb{F}=\{\mathcal{F}_t\}_{t\geq 0},P)$, 
then the Donsker delta function
$\delta_{B(t)}(x)$ of a Brownian motion $B$ at the point $x$ can be regarded as the time derivative of the local time $L_{t}(x)$ of $B$. More precisely, we have
    $$
L_{t}(x)=\int_{0}^{t}\delta_{B(s)}(x)ds.
    $$   
Such an integral exists as an element of the Hida space $(\mathcal{S}^{*})$ of stochastic distributions. See Section \ref{Hida}.

In \cite{LP} the authors use white noise theory to obtain an explicit solution formula for a general stochastic differential equation (SDE), and this is used
to find an expression for the Donsker delta function for the solution of an SDE.
Subsequently this was also extended to SDEs driven by L\'evy noise in \cite{MP}.

The main result of the current paper is that the Donsker delta measure of a McKean-Vlasov process (see below) always satisfies a certain Fokker-Planck type SPDE in the sense of distributions. Moreover, we use this to find explicit formulae for the Donsker delta functions for McKean-Vlasov processes, and hence their local times, in specific cases.

Let $X(t)=X_t \in \RR$ be the solution of a McKean-Vlasov SDE, i.e. a mean-field stochastic differential equation, of the form (using matrix notation),
\small
\begin{align*}
dX(t)=\alpha(t,X(t),\mu_t)dt+\beta(t,X(t),\mu_{t})dB(t); \quad
X(0) =Z\in\mathbb{R}.
\end{align*}
We call $X$ a \emph{McKean-Vlason process}.\\
Here the  $\sigma$-algebra $\mathbb{F}=\{\mathcal{F}_t\}_{t\geq0}$ denotes the filtration generated by $Z$ and $B(\cdot)$, $Z$ is a random variable which is independent of the $\sigma$-algebra generated by $B(\cdot)$ and such that $\EE[|Z|^2]<\infty.$\\
We denote by $\mu_t=\mu_t(\omega)=\mu_t(dx,\omega) $ the \emph{conditional law} of $X(t)$ given the filtration $\mathcal{F}_t$ generated by the Brownian motion $B$. More precisely, we consider the following model:

\begin{definition}
Define $\mu_t=\mu_{X(t)}(\omega, dx)$ to be regular conditional distribution of $X(t)$ given $\mathcal{F}_t$.
This means that $\mu_t(\omega,dx)$ is a Borel probability measure on $\RR$ for all $t \in [0,T],\quad \omega \in \Omega$ and
\begin{equation*} 
\int_{\mathbb{R}} g(x) \mu_t(dx,\omega)= \EE[g(X(t)) | \mathcal{F}_t](\omega),
\end{equation*}
 for all functions $g$ such that $\EE[ |g(X(t)) |] < \infty$.
\end{definition}


Since we consider only a one-dimensional Brownian motion $B(t) \in \RR$, we will show that the regular conditional distribution of $X(t)$ given the filtration $\{\mathcal{F}_t\}_{t\geq 0}$   can be identified with the Donsker delta measure in the sense of distribution. See details in Section \ref{sec3.1}

\section{Preliminaries}
In this section we review some basic notions and results that will be used throughout this work.
\subsection{Radon measures}
 A Radon measure on $\mathbb{R}^d$ is a Borel measure which is finite on compact sets, outer regular on all Borel sets and inner regular on all open sets. In particular, all Borel probability measures on $\mathbb{R}^d$ are Radon measures.\\
In the following, we let 
\begin{itemize}
    \item $\mathbb{M}_0$ be the set of deterministic Radon measures.   
    \item $C_0(\mathbb{R}^d)$ be the uniform closure of the space $C_c(\mathbb{R}^d)$ of continuous functions with compact support.
\end{itemize}
If we equip $\mathbb{M}_0$ with the total variation norm $||\mu||:=|\mu|(\mathbb{R}^d)$, then $\mathbb{M}_0$ becomes a Banach space, and it is the dual of $C_0(\mathbb{R}^d)$. See  Chapter 7 in Folland \cite{F} for more information.\\
If $\mu \in \mathbb{M}_0$ is a finite measure, we define
\begin{equation}\label{mu^}
\widehat{\mu}(y):=F[\mu](y):=\int_{\mathbb{R}^d}e^{-ixy} \mu(dx); \quad y \in \mathbb{R}^d
\end{equation}
to be the Fourier transform of $\mu $ at $y$. \\
In particular, if $\mu(dx)$ is absolutely continuous with respect to Lebesgue measure $dx$ with Radon-Nikodym-derivative $m(x)=\frac{\mu(dx)}{dx}$, so that $\mu(dx)=m(x) dx$ with $m \in L^1(\RR^d)$, we define
the Fourier transform
of $m$ at $y$, denoted by $\widehat{m}(y)$ or $F[m](y)$, by
\begin{align*}
F[m](y)=\widehat{m}(y)=\int_{\mathbb{R}^d} e^{-ixy} m(x) dx; \quad y \in \mathbb{R}^d.
\end{align*}

We let $\mathbb{M}$ denote the set of all random measures $\mu(dx,\omega); \omega \in \Omega$ such that $\mu(dx,\omega) \in \mathbb{M}_0$ for each given $\omega \in \Omega$.

\subsection{The Schwartz space of tempered distributions}
We recall now some notions from white noise analysis.
\begin{itemize}
    \item  $\mathcal{S}=\mathcal{S}(\mathbb{R}^d)$ be the
Schwartz space of rapidly decreasing smooth real
functions on $\mathbb{R}^d$. It is a
Fr\'echet space with respect to the family of seminorms:\label{simb-029} 
\begin{equation*}
\Vert f \Vert_{k,\alpha} := \sup_{x \in \mathbb{R}^d}\big\{ (1+|x|^k) \vert
\partial^\alpha f(x)\vert \big\},
\end{equation*}
where $k = 0,1,...$, $\alpha=(\alpha_1,...,\alpha_d)$ is a multi-index with $%
\alpha_j= 0,1,...$ $(j=1,...,d)$ and\label{simb-030} 
\begin{equation*}
\partial^\alpha f := \frac{\partial^{|\alpha|} f}{\partial
x_1^{\alpha_1}\cdots \partial x_d^{\alpha_d}}\quad \text{ for } |\alpha|=\alpha_1+ ... +\alpha_d.
\end{equation*}
    \item  $\mathcal{S}^{\prime }=\mathcal{S}^{\prime }(\mathbb{R}^{d})$ is the space of tempered distributions. It is  the dual of $\mathcal{S}$. 
\end{itemize}
\subsection{The Hida space $(\mathcal{S})^{*}$ of stochastic distributions}\label{Hida}
We restrict ourselves to the white noise probability space $(\Omega=\mathcal{S}^{\prime }, \mathcal{F}=\mathcal{B}, P)$, where $\mathcal{B}$ is the Borel $\sigma$-algebra and the probability $P$ is the probability measure on $\mathcal{S}^{\prime}$ defined in virtue of the Bochner-Minlos-Sazonov theorem).\\
Let $\mathcal{J}$ denote the set of all finite multi-indices 
$\alpha =(\alpha _{1},\alpha _{2},\ldots ,\alpha _{m})$, $m=1,2,\ldots$, of
non-negative integers $\alpha _{i}$. 
\begin{equation}  
\big( 2\mathbb{N} \big)^{\alpha} = \prod\limits_{j=1}^{m }(2j)^{\alpha
_{j}}=(2\cdot1)^{\alpha_1}(2\cdot 2)^{\alpha _2} (2\cdot 3)^{\alpha _3} ...
(2m)^{\alpha_m}.
\end{equation}
If $\alpha=(\alpha_1, \alpha_2, ...) \in \mathcal{J}$ we put \label{simb-036} 
\begin{equation}  \label{3.12a}
H_{\alpha }(\omega ):=\prod\limits_{j=1}^{m}h_{\alpha _{j}}(\theta
_{j}(\omega)) = h_{\alpha_1}(\theta_1) h_{\alpha_2}(\theta_2) ...
h_{\alpha_m}(\theta_m), \quad \omega\in\Omega.
\end{equation}

The family $\{H_{\alpha }\}_{\alpha \in 
\mathcal{J}%
\text{ }}$ constitutes an orthogonal basis of $L^{2}(P)$.

\begin{itemize}
    \item $((\mathcal{S})_{k})_{k\in\mathbb{R}}$ is the Hilbert space consisting  of all $f=\sum\limits_{\alpha
\in \mathcal{J}}c_{\alpha } H_{\alpha }\in L^{2}(P )$ such that 
$\Vert f \Vert^2_k := \sum\limits_{\alpha \in \mathcal{J}}\alpha!c_{\alpha
}^{2}(2\mathbb{N})^{\alpha k}<\infty,$ for numbers $c_{\alpha }\in \mathbb{R}$.

    \item The space $(\mathcal{S})= \bigcap_{k\in\mathbb{R}} (\mathcal{S})_{k}$ equipped with the projective topology is the Hida space of stochastic test functions.
\end{itemize}

\begin{itemize}
    \item $((\mathcal{S})_{-k})_{k\in\mathbb{R}}$ is the Hilbert space consisting of all formal sums $F=\sum\limits_{%
\alpha \in \mathcal{J}} c_{\alpha }H_{\alpha }$  equipped with the norm
\begin{equation*}
\Vert F \Vert^2_{-k} :=\sum\limits_{\alpha \in \mathcal{J}}\alpha
!c_{\alpha}^{2} (2\mathbb{N})^{-\alpha k}<\infty.
\end{equation*}
    \item The space $(\mathcal{S})^{*} = \bigcup_{k\in\mathbb{R}} (\mathcal{S})_{-k}$ equipped with the inductive topology is the Hida space of stochastic distributions. It can be regarded as the dual of $(\mathcal{S})$.
\end{itemize}

\subsection{The Donsker delta function}
We now recall some basic definitions:
\begin{definition}
Let $Y :\Omega\rightarrow\mathbb{R}$ be a random variable which also belongs to the Hida space $(\mathcal{S})^{\ast}$ of stochastic distributions. Then a continuous function
\begin{equation}\label{donsker}
    \delta_Y(\cdot): \mathbb{R}\rightarrow (\mathcal{S})^{\ast}
\end{equation}
is called a Donsker delta function of $Y$ if it has the property that
\begin{equation}\label{donsker property }
    \int_{\mathbb{R}}g(y)\delta_Y(y)dy = g(Y) \quad a.s.
\end{equation}
for all (measurable) $g : \mathbb{R} \rightarrow \mathbb{R}$ such that the integral converges in $(\mathcal{S})^{*}.$
\end{definition}
The Donsker delta function is related to the \emph{regular conditional distribution}. The connection is the following:
The \emph{regular conditional distribution} with respect to the $\sigma$-algebra $\mathcal{F}
$ of a given real random variable $Y$, denoted by $\mu_Y(dy)=\mu_Y(dy, \omega); \omega \in \Omega$, is defined by the following properties:
\begin{itemize}
\item
For any Borel set $\Lambda \subseteq \mathbb{R}$,  $\mu_Y(\Lambda,\cdot)$ is a version of $\mathbb{E}[\mathbf{1}_{Y \in \Lambda} | \mathcal{F}]$.
\item
For each fixed $\omega \in \Omega$, $\mu_Y(dy, \omega)$ is a probability measure on the Borel subsets of $\mathbb{R}$.
\end{itemize}

It is well-known that such a regular conditional distribution always exists. See e.g. \cite{B}, p.79.\\
From the required properties of $\mu_Y(dy, \omega)$, we get the following formula:
\begin{equation}
\int_{\mathbb{R}} f(y) \mu_Y(dy,\omega) = \mathbb{E}[ f(Y) | \mathcal{F}].
\end{equation}

\begin{definition}
We call $\mu_Y(dy,\omega)$ the \textbf{Donsker delta measure} of the random variable $Y$ and denote it by $\delta_{Y}(dy,\omega)$.
\end{definition}
Comparing this with the definition of the Donsker delta function, we obtain the following representation of the regular conditional distribution:

\begin{lemma}
Suppose $\mu_Y(dy, \omega)$ is absolutely continuous with respect to Lebesgue measure $dy$ on $\mathbb{R}$ and that $Y$ is measurable  with resepct to $\mathcal{F}$. Then the \textbf{Donsker delta function} of $Y$,  $\delta_Y(y,\omega),$ is the Radon-Nikodym derivative of $\mu_{Y}(dy,\omega)$ with respect to Lebesgue measure $dy$, i.e.
\begin{equation}
\delta_Y(y,\omega)=\frac{\mu_Y(dy, \omega)}{dy}.
\end{equation}
\end{lemma}

We will prove in Theorem \ref{DonskerFP} that the Donsker delta function can be regarded as a stochastic distribution in $\mathcal{S}^{'}$, satisfying a Fokker-Planck type SPDE in the sense of distributions. It can also be represented as an element of the Hida stochastic distribution space $(\mathcal{S})^{*}$, and as such it can in some cases be expressed explicitly in terms of Wick calculus. For example, if $Y(t)=B(t)$, we have
\begin{align}\label{donsker3}
\delta_{B(t)}(x)= (2 \pi t)^{-\tfrac{1}{2}} \exp^{\diamond}\Big(-\frac{(B(t)-x)^{\diamond 2}}{2t}\Big) \in (\mathcal{S})^{*},
\end{align}
where $\diamond$ denotes Wick multiplication and $\exp^{\diamond}$ denotes Wick exponential.
Note that even though the Donsker delta function can only be represented as a distribution, its conditional expectation  can be a real valued stochastic process. For example, for $t < T$  we have
\begin{align}
  \mathbb{E}[\delta_{B(T)}(x)|\mathcal{F}_t]
   &= (2\pi (T-t))^{-\frac{1}{2}} \exp\Big[- \frac{(B(t)-x)^2}{2(T-t)}\Big].
\end{align}
For more examples, we refer to e.g. \cite{AOU} or \cite{DOP}.

\section{The Donsker delta equation for McKean-Vlasov processes}
\subsection{ The general multidimensional Fokker-Planck equation}\label{sec3.1}
To explain the background for this section, let us recall the general multidimensional situation studied in \cite{AO}, where $X(t) \in \mathbb{R}^d$ is a McKean-Vlasov diffusion, of the form (using matrix notation),
\begin{align}\label{MV}
dX(t)& =b(t,X(t),\mu_t)dt+\sigma(t,X(t),\mu_{t})dB(t),\quad
X(0) =Z,
\end{align}
where $B $ is a multi-dimensional Brownian motion. \\
Here $Z$ is a random variable which is independent of the $\sigma$-algebra generated by $B(\cdot)$ and such that
$$
\EE[|Z|^2]<\infty.
$$
Define the $\sigma$-algebra $\mathcal{F}=\{\mathcal{F}_t\}_{t\geq0}$ to be the filtration generated by $Z$ and $B(\cdot)$.
\vskip 0.2cm
Let $\mathbb{M}$ denote the set of all Borel measures on $\mathbb{R}^d$. We assume that the coefficients $b(t,x,\mu):[0,T]\times \mathbb{R}^d \times \mathbb{M}\rightarrow \mathbb{R}^d$ and $\sigma(t,x,\mu):[0,T]\times \mathbb{R}^d \times \mathbb{M}\rightarrow \mathbb{R}^d$  are bounded and $\mathbb{F}$-predictable processes for all $x,\mu$, 
and that $b$ and $ \sigma$ are continuous with respect to $t$ and $x$ for all $\mu$. 
\vskip 0.2cm
One can check that under some assumptions, such as Lipschitz and linear growth conditions, there exists a unique solution of equation  \eqref{MV}. 
\vskip 0.2cm
We denote by $\mu_t=\mu_t(\omega)=\mu_t(dx,\omega) $ the \emph{conditional law} of $X(t)$ given the filtration $\mathcal{F}_t$ generated by the Brownian motion $B$. More precisely, we consider the following model:

\begin{definition}\label{Donsker 1}
Fix one of the Brownian motions, say $B_1=B_1(t,\omega)$, with filtration $\{\mathcal{F}_t^{(1)}\}_{t\geq 0}$. We define $\mu_t=\mu_{X(t)}(\omega, dx)$ to be regular conditional distribution of $X(t)$ given $\mathcal{F}_t^{(1)}$.
This means that $\mu_t(\omega,dx) $ is a Borel probability measure on $\RR^d$ for all $t \in [0,T],\omega \in \Omega$ and
\begin{equation} \label{cond}
\int_{\mathbb{R}^d} g(x) \mu_t(dx,\omega)= \EE[g(X(t)) | \mathcal{F}_t^{(1)}](\omega)
\end{equation}
 for all functions $g$ such that $\EE[ |g(X(t)) |] < \infty$.
\end{definition}

The following version of the stochastic Fokker-Planck integro-differential equation for the conditional law  for McKean-Vlasov jump diffusions was proved by Agram and \O ksendal \cite{AO}. For simplicity we consider only the case without jumps here.

\begin{theorem}{(Conditional stochastic Fokker-Planck equation \cite{AO})}\\
Let $X(t)$ be as in \eqref{MV} with $m\geq 2$ and let $\mu_t:=\mu_{X(t)}(dx,\omega)$ be the regular conditional distribution of $X(t)$ given $\mathcal{F}_t^{(1)}$.

Then for a.a. $\omega \in \Omega$ the conditional law $\mu_t \in \mathcal{S}^{'}$ and it satisfies the following SPDE (in the sense of distributions):

 \begin{align} 
d\mu _{t} =A_0^{*} \mu_t dt + A_1^{*}\mu_t dB_1(t),  \quad \mu_0=\mathcal{L}(X(0)).\label{FP}
\end{align}
Here 
$A_0^{*}, A_1^{*}$ are the integro-differential operator and the differential operator which are given respectively by:
\begin{align}
A_0^{*}\mu&= -\sum_{j=1}^d D_j [b_j \mu] +\frac{1}{2}\sum_{n,j=1}^d D_{n,j}[(\sigma \sigma^{T})_{n,j} \mu], \nonumber\\ \label{A0*}
\end{align}
and
\begin{align}
A_1^{*}\mu= - \sum_{j=1}^d D_j[\beta_{1,j} \mu], \label{A1*}
\end{align}

In the above $D_j, D_{n,j}$
denote $\frac{\partial }{\partial x_j}$  and  $\frac{\partial^2}{\partial x_n \partial x_j}$ respectively, in the sense of distributions.
\end{theorem}

\subsection{The Fokker-Planck equation for the Donsker measure}
In \cite{AO} the theorem above was proved under the assumption that $m \geq 2$. However, the proof also works if $m=1$ and $\mathcal{F}_t^{(1)}=\mathcal{F}_t$.
Note that in this case, since $X(t)$ is $\mathcal{F}_t$-measurable, the identity \eqref{cond} states that
\begin{equation} \label{cond2}
\int_{\mathbb{R}^d} g(x) \mu_t(dx,\omega)= g(X(t)) 
\end{equation}
 for all functions $g$ such that $\int_{\mathbb{R}^d} |g(x)| \mu_t(dx,\omega) < \infty$.

In particular, if we choose $d=m=1$ in the above we get that the conditional law coincides with the Donsker measure, i.e.
\begin{align}
\mu_t(x,\omega)= \delta_{X(t)}(dx,\omega).
\end{align}
Therefore we get the following Fokker-Planck equation for the Donsker measure:
\begin{theorem}\label{DonskerFP}
Assume that $X(t)$ is as in \eqref{MV}, but with $d=m=1$.\\
Then the Donsker delta measure $\mu_t= \delta_{X(t)}(dx,\omega)$ satisfies the following equation (in the sense of distribution)
 \begin{align} 
d\mu _{t} &=\Big\{-D[b(t,x,\mu_t) \mu_t] +\tfrac{1}{2}D^2[\sigma^2(t,x,\mu_t) \mu_t] \Big\} dt - D[\sigma(t,x,\mu_t) \mu_t] dB(t);\quad t\geq 0\nonumber\\
  \mu_0&=\mathcal{L}(X(0)),\label{A}
\end{align}
where $D=\frac{\partial}{\partial x}$ and $D^{2}=\frac{\partial ^{2}}{\partial x^{2}}$.
\end{theorem}

\section{Local time}

In this section we first recall the definition of local time of a stochastic process $Y(\cdot)$:

\begin{definition}
The local time $L_{t}(y)$ of $Y(\cdot)$ at the point $y$ and at time
$t$ is defined by
\[
L_{t}(y)=\lim_{\epsilon\rightarrow0}\frac{1}{2\epsilon}\lambda(\{s\in
\lbrack0,t];Y(s)\in(y-\epsilon,y+\epsilon)\}),
\]
where $\lambda$ denotes Lebesgue measure on $\mathbb{R}$ and the limit is in $L^2(P)$.
\end{definition}


In the white noise context the local time can be represented as the integral of the Donsker delta function. More precisely, we have the following result:
\begin{theorem}
The local time $L_{t}(x)$ of $X$ at the point $x$ and the time $t$
is given by
\begin{equation}
L_{t}(z)=\int_{0}^{t}\delta_{X(s)}(x)ds, \label{eq2.7}%
\end{equation}
where the integration takes place in $(\mathcal{S})^{*}$ (or in $\mathcal{S}'$ for each $\omega$).
\end{theorem}

\dproof For completeness we give the proof.\\By definition of the local time and the Donsker delta function, we have
\begin{align*}
L_{t}(z)  &  =\lim_{\epsilon\rightarrow0} \int_{0}^{t} \chi_{(z-\epsilon
,z+\epsilon)} (Y(s))ds \\
&=\lim_{\epsilon\rightarrow0} \int_{0}^{t}
\Big(\int_{\mathbb{R}}\chi_{(z-\epsilon,z+\epsilon)} (y) \delta_{Y(s)}(y)
dy\Big)ds\\
&  =\lim_{\epsilon\rightarrow0} \int_{\mathbb{R}}\chi_{(z-\epsilon
,z+\epsilon)} (y) \Big(\int_{0}^{t} \delta_{Y(s)}(y) ds\Big)dy =\int_{0}^{t}
\delta_{Y(s)}(z) ds,
\end{align*}
because the function $y \mapsto\delta_{Y(s)}(y)$ is continuous in
$(\mathcal{S})^{*}$ (and in $\mathcal{S}'$).
\fproof

\begin{remark}
Note that even though we in general can only say that $\delta_{X(t)}(x) \in (\mathcal{S})^{*}$,  $L_t(x)$ usually exists as a real-valued stochastic process.
\end{remark}

\section{Explicit solutions}
In this Section, we find explicitly the Donsker delta function for some particular McKean-Vlasov processes and accordingly their local time.\\
Suppose that $\mu _{t}$ is absolutely continuous i.e.
\begin{equation}
\mu _{t}(dx)=m(t,x)dx. \label{abs}
\end{equation}
Then \eqref{MV} gets the form
\begin{align}
dX(t)& =b(t,X(t),m_t)dt+\sigma(t,X(t),m_{t})dB(t);\quad
X(0) =Z,
\end{align}
where $m_t(x)= m(t,x)$
and 
\eqref{A} becomes a \emph{stochastic partial differential equation (SPDE)}, as follows:
\begin{theorem}
Suppose \eqref{abs} holds. Then the Donsker delta function $m(t,x)=\delta_{X(t)}(x)$ is the solution in $(\mathcal{S})^{*}$ of the following SPDE:
 \begin{align} 
d_t m(t,x) &=\Big\{-\tfrac{\partial}{\partial x}[b(t,x,m) m(t,x)] +\tfrac{1}{2}\tfrac{\partial ^2}{\partial x ^2} [\sigma^2(t,x,m) m(t,x)] \Big\} dt  \label{C}\\
&- \tfrac{\partial}{\partial x}[\sigma(t,x,m) m(t,x)] dB(t);\quad t\geq 0,\nonumber\\
  m(0,x)&=\tfrac{\partial}{\partial x} \mathcal{L}(X(0)).
\end{align}
\end{theorem}
\subsection{Brownian motion}
Consider the special case when $X(t)=B(t); B(0)=Z$. Then $b=0$ and $\sigma=1$ and  equation \eqref{A} becomes
\begin{align}  
\frac{\partial m}{\partial t}(t,x)&=\frac{1}{2}\frac{\partial^{2}m}{\partial x^{2}
}(t,x)+\frac{\partial m}{\partial x}(t,x)\diamond \dot{B}(t);\quad t\geq 0,\label{B}\\
  m(0,x)&=\tfrac{\partial}{\partial x} \mathcal{L}(X(0))\label{B0}.
\end{align}
We can easily verify by Wick calculus that a solution in $(\mathcal{S})^{*}$ of equation \eqref{B}  is
\begin{align}\label{donsker4}
\delta_{B(t)}(x)= (2 \pi t)^{-\tfrac{1}{2}} \exp^{\diamond}\Big(-\frac{(B(t)-x)^{\diamond 2}}{2t}\Big),
\end{align}
which is in agreement with \eqref{donsker3}. The details are as follows:\\
Try
\[
m(t,x)=\frac{1}{\sqrt{2\pi t}}\exp^{\diamond}[-\frac{(x-B(t))^{\diamond2}}{2t}].
\]
Then%
\begin{align*}
\frac{\partial m}{\partial t} (t,x) & =-\frac{1}{2}\frac{t^{-3/2}}{\sqrt{2\pi}}%
\exp^{\diamond}[-\frac{(x-B(t))^{\diamond2}}{2t}]+\frac{1}{\sqrt{2\pi t}}\exp^{\diamond}[-\frac{(x-B(t))^{\diamond2}}{2t}](-\frac{x-B(t)}%
{ t})\diamond\dot{B}(t)\\
& +\frac{1}{\sqrt{2\pi t}}\exp^{\diamond}[-\frac{(x-B(t))^{\diamond2}}{2t}]\frac{(x-B(t))^{2}}{2t^{2}},
\end{align*}
and%
\[
\frac{\partial m}{\partial x}=\frac{1}{\sqrt{2\pi t}}\exp^{\diamond
}[-\frac{(x-B(t))^{\diamond2}}{2t}]\diamond(-\frac{x-B(t)}{t}),
\]
and%
\[
\frac{\partial^{2}m}{\partial x^{2}}(t,x)=\frac{1}{\sqrt{2\pi t}}\exp^{\diamond
}[-\frac{(x-B(t))^{\diamond2}}{2t}]\diamond(\frac{x-B(t)}{t})^{\diamond2}+\frac{1}{\sqrt{2\pi t}}%
\exp^{\diamond}[-\frac{(x-B(t))^{\diamond2}}{2t}](-\frac{1}{t}).
\]
Collecting the terms we see that%
\[
m(t,x)=\frac{1}{\sqrt{2\pi t}}\exp^{\diamond}(-\frac{(x-B(t))^{\diamond2}}%
{2t}),
\]
satisfies the Fokker-Planck equation \eqref{B} for the conditional law of $B(t)$.

From white noise theory we know that
\begin{itemize}
\item
$E[X \diamond Y]=E[X] E[Y]$
\item
$E[\exp^{\diamond}Y]:=E[\sum_{n=0}^{\infty} \tfrac{1}{n!} Y^{\diamond n}]=\sum_{n=0}^{\infty} \tfrac{1}{n!} E[Y^{\diamond n}]= \sum_{n=0}^{\infty} \tfrac{1}{n!} E[Y]^{ n}]=\exp( E[Y])$
\end{itemize}
for all random variables $X,Y$ with a finite expectation (independent or not). From this we see that
\begin{align}
    E[\delta_{B(t)}(x)]=\frac{1}{\sqrt{2\pi t}} \exp \big(E[-\frac{(x-B(t))^{\diamond 2}}{2t}]\big)=\frac{1}{\sqrt{2\pi t}} \exp\Big( -\frac{(x-E[Z])^2}{2t}\big).
\end{align}
In particular, if $X(0)=Z=z$ (constant) $\in \mathbb{R}$ a.e., then
\begin{align} \label{E}
    E[\delta_{B(t)}(x)]=\frac{1}{\sqrt{2\pi t}} \exp\Big( -\frac{(x-z)^2}{2t}\big),
\end{align}
which has a singularity at $x=z$.
\subsection{Coefficients not depending on $x$}
The next result shows that, under some conditions, the Donsker delta function can be an ordinary function if the initial value $X(0)$ has a density: 
\begin{theorem}\label{4.2}
Assume that $X(t)$ is the solution of the following McKean-Vlasov equation:
\begin{align}\label{MV1}
dX(t)=\alpha(t,\mu_t)dt+\beta(t,\mu_{t})dB(t); \quad
X(0) =Z,
\end{align}
where the coefficients $\alpha(t,x,\mu)=\alpha(t,\mu)$ and $\beta(t,x,\mu)=\beta(t,\mu)$ do not depend on $x$.
Suppose that $X(0)=Z$ is a random variable (independent of $B$) with density 
\begin{align}
h(z)=\frac{\partial }{\partial z}\mathcal{L}(Z)(z); \quad  z \in \mathbb{R}.
\end{align}
\begin{enumerate}
\item
Define
\begin{align}\label{m}
Y_t(x)=h\big(x - \int_0^t \alpha(s,\mu_s)ds - \int_0^t \beta(s,\mu_s) dB(s)\big). \\
\end{align}
Then $Y(t,x)$ is the Donsker delta function of $X(t)$.
\item
The solution $X(t)$ of \eqref{MV1} is given by
\begin{align}
X(t)= \int_{\mathbb{R}} x Y_t(x)dx.
\end{align}
\end{enumerate}
\end{theorem}

\dproof
\begin{enumerate}
    \item We show that $Y(t,x)$ satisfies equation \eqref{C}.\\
By the Ito formula we have
\begin{align}
d_t Y(t,x)&= h'(Z(t,x)) d_tZ(t) +\tfrac{1}{2} h''(Z(t,x)) \beta^2(t) dt\nonumber\\
&=\{-\alpha(t,\mu_t) h'(Z(t,x)) + \tfrac{1}{2} \beta^2(t,\mu_t) h''(Z(t,x)) \} dt -\beta(t,\mu_t) h'(Z(t,x)) dB(t). \label{dY}
\end{align}
Since
\begin{align}
h'(Z(t,x))=\tfrac{d}{dz}h(z)_{z=Z(t,x)}= \tfrac{\partial}{\partial x} Y(t,x),
\end{align}
we see that the equation \eqref{dY} can be written
\begin{align}
d_t Y(t,x)=[-\alpha(t,\mu_t) \tfrac{\partial}{\partial x} Y(t,x) + \tfrac{1}{2} \beta^2(t,\mu_t) \tfrac{\partial^2}{\partial x^2} Y(t,x) ] dt -\beta(t,\mu_t) \tfrac{\partial}{\partial x} Y(t,x) dB(t),
\end{align}
which is the same as equation \eqref{C}.\\
Since $Y(0,x)=h(x)=m(0,x)$ we conclude by uniqueness that $Y(t,x)=m(t,x)$ for all $t$.\\
\item This follows from the definition of the Donsker delta function.
\end{enumerate}
\fproof

\subsubsection{Constant coefficients}
As a special case of the case above, suppose that
\begin{align}\label{MK}
dX(t)& =\alpha dt+\beta dB(t),\quad 
X(0) =Z,
\end{align}
where $\alpha$ and $\beta$ are constants. Then by Theorem \ref{4.2} the Donsker delta function is 
\begin{align}
\delta_{X(t)}(x)&=h(x -\alpha t -\beta B(t)).
\end{align}

\subsection{Mean-field geometric Brownian motion}
Suppose that $X(t)$ is a McKean-Vlasov process of the form
\begin{equation}
dX_{t}=\alpha(t,\mu_t)X_{t}dt+\beta(t,\mu_t) X_{t}dB_{t};\quad X_{0}=Z > 0.  \label{gbm}
\end{equation}%
We call this a \emph{mean-field geometric Brownian motion}. For such processes we have:
\begin{theorem}
\begin{myenumerate}
\item
 The Donsker delta function $m_t(x)$ for the mean-field geometric Brownian motion $X(t)$ is
\begin{align}
m_t(x)=\delta _{X_{t}}(x)=\frac{1}{x} H\big(\ln{x} - \int_0^t \alpha(s,\mu_s)ds - \int_0^t \beta(s,\mu_s) dB(s)\big),
\end{align}
where
\begin{align}
H(z)=\frac{\partial}{\partial z} \mathcal{L}(\ln{Z})(z); \quad z \in \mathbb{R}.
\end{align}
\item
 The solution $X(t)$ of the mean-field geometric Brownian motion equation \eqref{gbm} can be written
\begin{align}
X(t)&=\int_0^\infty H\big(\ln{x} - \int_0^t \alpha(s,\mu_s)ds - \int_0^t \beta(s,\mu_s) dB(s)\big)dx\nonumber\\
&= \int_{\mathbb{R}} e^{u}H\big(u - \int_0^t \alpha(s,\mu_s)ds - \int_0^t \beta(s,\mu_s) dB(s)\big)du.
\end{align}
\end{myenumerate}
\end{theorem}

\dproof \begin{itemize}
    \item [(i) ] The corresponding Fokker-Planck equation for
the Donsker delta function\\ $m_{t}(x)=\delta _{X(t)}(x)$ is
\small
\begin{align}
dm_{t}(x) &=\{-\frac{\partial }{\partial x}[\alpha(t,m_t)xm_{t}(x)]+\frac{1}{2}\beta^2(t,m_t)
\frac{\partial ^2 }{\partial x ^2}[x^{2}m_{t}(x)]\}dt-\beta(t,m_t) \frac{\partial }{\partial x}[xm_{t}(x)]dB_{t}  \nonumber \\
&=\{-\alpha(t,m_t) m_{t}(x)-\alpha(t,m_t) xm_{t}^{\prime }(x)+\frac{1}{2}\beta^2(t,m_t)
[2m_{t}(x)+4xm_{t}^{\prime }(x)+x^{2}m^{\prime \prime }(x)]\}dt  \nonumber\\
&-\beta(t,m_t) \lbrack m_{t}(x)+xm_{t}(x)]dB_{t}; \quad m_0(x)=\frac{\partial}{\partial x}\mathcal{L}(Z)(x).  \label{4.9}
\end{align}%
This is a stochastic partial differential equation in $m_{t}(x).$ It seems difficult to find directly an explicit solution of this equation.
However, we can find the solution $m_t(x)=\delta _{X(t)}(x)$ by proceeding as follows:\\
The solution of (\ref{gbm}) is%
\[
X_{t}=Z \exp \big(\int_0^t\beta(s,m_s) dB(s)+\int_0^t \{\alpha(s,m_s) -\frac{1}{2}\beta ^{2}(s,m_s)\}ds\big)=\exp (Y_{t}),
\]%
where%
\[
Y_{t}=\ln{Z}+\int_0^t\beta(s,m_s) dB(s)+\int_0^t \{\alpha(s,m_s) -\frac{1}{2}\beta ^{2}(s,m_s)\}ds.
\]
By Theorem \ref{4.2} we know that%
\[
\delta _{Y_{t}}(x)=H\big(x - \int_0^t \alpha(s,\mu_s)ds - \int_0^t \beta(s,\mu_s) dB(s)\big),
\]
where%
\[
H(z)=\frac{\partial}{\partial z} \mathcal{L}(\ln{Z})(z).
\]%
By definition \ we have%
\[
\int_{\mathbb{R}^{+}}g(y)\delta _{Y_{t}}(y)dy=g(Y_{t}).
\]%
With $g(y)=\exp (y)$ this gives
\[
\int_{\mathbb{R}}\exp(y)\delta _{Y_{t}}(y)dy=\exp (Y_{t})=X_{t}.
\]
Hence, substituting $\exp(y)=x$,%
\[
X_{t}=\int_{\mathbb{R}}\exp(y)\delta _{Y_{t}}(y)dy=\int_{\mathbb{R}^{+}}
x\delta _{Y_{t}}(\ln(x))\frac{dx}{x}.
\]%
From this we deduce that
\begin{align}
m_t(x)=\delta _{X_{t}}(x)=\frac{\delta _{Y_{t}}(\ln(x))}{x}=\frac{1}{x} H\big(\ln{x} - \int_0^t \alpha(s,\mu_s)ds - \int_0^t \beta(s,\mu_s) dB(s)\big)
\end{align}
is the Donsker delta function of $X_{t}$.\\
    \item [(ii) ]  This part follows by the definition of the Donsker delta function.
\end{itemize}
\fproof


\subsection{An example related to the Burgers equation}
Suppose the McKean-Vlasov equation has the form
\begin{align} \label{Burger}
dX(t)&=\alpha m(t,X(t))dt + \beta dB(t); \quad
X(0) = Z,
\end{align}
where $m(t,x)=\frac{\partial }{\partial x}\mu(t,x)=\frac{\partial }{\partial x}\mathcal{L}(X(t))(x)\in L^2([0,T] \times \mathbb{R}).$\\
Then the corresponding FP equation for the Donsker function $m(t,x)$ is
\begin{align}
dm(t,x)&=\{ -\alpha \frac{\partial}{\partial x}(m^2(t,x)) + \frac{1}{2} \beta^2 \frac{\partial ^2}{\partial  x^2} m(t,x)\} dt - \beta \frac{\partial }{\partial x}m(t,x) dB(t), \label{4.7}\\
m(0,x)&=h(x)= \frac{\partial}{\partial x}Z(x). 
\end{align}
This is a stochastic Burgers equation. It is well-known that by using the Cole-Hopf transformation the equation can be transformed into the classical heat equation. The details are as follows:
If we introduce a new function $\psi=\psi(t,x)$ such that
\begin{align}
m:= \psi_x := \frac{\partial}{\partial x} \psi, 
\end{align}
then we see that 
the Burgers equation \eqref{4.7} becomes the following equation in $\psi$:
\begin{align}
(\psi_x)_t=-2\alpha \psi_x (\psi_x)_x + \frac{1}{2} \beta^2(\psi_x)_{xx}-\beta \psi_{xx} \diamond \dot{B}(t).
\end{align}
Integrating with respect to $x$ this gives
\begin{align}
\psi_t=-\alpha (\psi_x)^2 + \frac{1}{2}\beta^2 \psi_{xx}-\beta \psi_{x} \diamond \dot{B}(t).
\end{align}
Now define the function $\varphi=\varphi(t,x)$ by
\begin{equation}
\psi= \gamma \ln{\varphi},
\end{equation}
for some constant $\gamma$.
Then in terms of $\varphi$ the above equation gets the form
\begin{align}
\gamma \frac{\varphi_t}{\varphi}&= -\alpha \gamma^2(\frac{\varphi_x}{\varphi})^2 +\frac{1}{2} \gamma \beta^2(\frac{\varphi_x}{\varphi})_x -\beta \gamma\frac{\varphi_x}{\varphi}\diamond \dot{B}(t)\nonumber\\
&=  -\alpha \gamma^2 (\frac{\varphi_x}{\varphi})^2 + \tfrac{1}{2} \gamma \beta^2 (\frac{\varphi \varphi_{xx}- (\varphi_x)^2}{\varphi^2})_x -\beta \gamma \frac{\varphi_x}{\varphi}\diamond \dot{B}(t)\nonumber\\
&=  -\alpha \gamma^2 (\frac{\varphi_x}{\varphi})^2 + \frac{1}{2} \gamma \beta^2 (\frac{\varphi_{xx}}{\varphi}) - \frac{1}{2}\gamma \beta^2(\frac{\varphi_x}{\varphi})^2 -\beta \gamma \frac{\varphi_x}{\varphi}\diamond \dot{B}(t).
\end{align}
This simplifies to
\begin{align}
\varphi_t= -(\gamma \alpha+\frac{1}{2}\beta^2)(\frac{\varphi_x}{\varphi})^2+\frac{1}{2} \beta^2 \varphi_{xx}-\beta \varphi_x \diamond \dot{B}(t).
\end{align}
If we choose
\begin{equation}
\gamma=-\frac{\beta^2}{2\alpha},
\end{equation}
the equation for $\varphi$ reduces to the (linear) stochastic heat equation
\begin{align}
\varphi_t&= \frac{1}{2} \beta^2 \varphi_{xx}-\beta \varphi_x \diamond \dot{B}(t),\label{4.16}\\
\varphi(0,x)&=k(x) \text{ (to be determined)},
\end{align}
or, using Ito differential notation,
\begin{align}
d\varphi(t,x)&= \frac{1}{2} \beta^2 \varphi_{xx}(t,x)dt-\beta \varphi_x(t,x) dB(t); \quad t 
\geq 0,\\
\varphi(0,x)&=k(x).
\end{align}
To find an expression for the solution of \eqref{4.16}, define an auxiliary process $R(t)=R^{(x)}(t)$ by
\begin{align}
R(t)=x - \beta B(t) + \beta \tilde{B}(t);\quad t \geq 0,
\end{align}
where $\tilde{B}$ is an auxiliary Brownian motion with law $\tilde{P}$ and independent of $B$.
Then by the Feynman-Kac formula
\begin{align}
\varphi(t,x):= \tilde{E}[k(R^{(x)} (t))]=\tilde{E}[k(x-\beta B(t) +\beta \tilde{B}(t))],
\end{align}
where $\tilde{E}$ denotes expectation with respect to $\tilde{P}$ and $k(z)=\varphi(0,z)$, solves equation \eqref{4.16}.
Going back to $m$ we get
\begin{align}
m(t,x)= \psi_x (t,x)= \gamma \frac{\partial}{\partial x} \ln{\varphi(t,x))}= \gamma\frac{\varphi_x (t,x)}{\varphi(t,x)}=\gamma\frac{\tilde{E}[k_x (x-\beta B(t) +\beta \tilde{B}(t))]}{\tilde{E}[k(x-\beta B(t) +\beta \tilde{B}(t))]}.
\end{align}
In particular, setting $t=0$ we get
\begin{align}
h(x):=m(0,x)= \gamma\frac{\tilde{E}[k_x (x)]}{\tilde{E}[k(x)]}= \gamma \frac {k_x(x)}{k(x)},
\end{align}
from which we deduce that
\begin{align}
k(x)= \exp(\frac{1}{\gamma} \int_0^x h(y)dy).
\end{align}
We summarize what we have proved as follows:

\begin{theorem}
\begin{enumerate}
    \item 
    The Donsker delta function $m(t,x)=\delta_{X(t)}(x)$ for the solution $X(t)$ of the McKean-Vlasov equation \eqref{Burger} is given by
    \begin{align}
m(t,x)= \gamma\frac{\varphi_x (t,x)}{\varphi(t,x)}=\gamma\frac{\tilde{E}[k_x (x-\beta B(t) +\beta \tilde{B}(t))]}{\tilde{E}[k(x-\beta B(t) +\beta \tilde{B}(t))]}, 
    \end{align}
    where
 \begin{align}
k(x)= \exp(\frac{1}{\gamma} \int_0^x m(0,y)dy)= \exp(\frac{1}{\gamma} \int_0^x \mathcal{L}(Z)(y)dy).
\end{align} 
and 
\begin{equation}
\gamma=-\frac{\beta^2}{2\alpha},
\end{equation}
\item
The solution $X(t)$ of \eqref{Burger} is given by
\begin{align}
X(t)=\int_{\mathbb{R}} x m(t,x) dx
\end{align}
with $m(t,x)$ as in part 1.
\end{enumerate}
\end{theorem}

\subsection{A solution approach based on Laplace and Fourier transforms}
Consider the Fokker-Planck equation
\begin{equation}
d\mu _{t}=\{-D[\alpha \mu ]+\frac{1}{2}D^{2}[\beta ^{2}\mu ]\}dt--D[\beta
\mu ]dB(t);\quad \mu _{0}=\delta _{x_{0}},  \label{(1)}
\end{equation}
for the McKean-Vlasov equation \eqref{MK}.
If $\alpha ,\beta $ are constants, this becomes
\begin{equation}
d\mu _{t}=\{-\alpha D[\mu ]+\frac{1}{2}\beta ^{2}D^{2}[\mu ]\}dt-\beta D[\mu
]dB(t);\quad \mu _{0}=\delta _{x_{0}}.  \label{(2)}
\end{equation}
If $d\mu _{t}=m(t,x)dx$, the equation can be written as
\begin{equation}
\frac{\partial }{\partial t}m(t,x)=-\alpha \frac{\partial }{\partial x}%
m(t,x)+\frac{1}{2}\beta ^{2}\frac{\partial ^{2}}{\partial x^{2}}m(t,x)+\beta 
\frac{\partial }{\partial x}m(t,x)\diamond \mathring{B}(t).  \label{(3)}
\end{equation}%
Let%
\begin{equation}
\widetilde{f}(s)=Lf(s)=\int_{0}^{s}e^{-st}f(t)dt\text{ \ denote the Laplace
transform},  \label{(4)}
\end{equation}
and
\begin{equation}
\widehat{f}(y)=Ff(y)=\int_{\mathbb{R}}e^{-ixy}f(x)dx\text{ \ denote the
Fourier transform}.  \label{(5)}
\end{equation}%
Then 
\begin{equation}
L\Big(\frac{\partial }{\partial t}f(t)\Big)(s)=s(Lf)(s)-f(0),  \label{(6)}
\end{equation}
and%
\begin{equation}
L(\exp (bt))(s)=\frac{1}{s-b},  \label{(7)}
\end{equation}%
and%
\begin{equation}
F[D^{n}w](y)=(iy)^{n}F[w](y).  \label{(8)}
\end{equation}%
Hence, applying the Laplace and Fourier transform to (\ref{(3)}), we get
\small
\begin{eqnarray*}
s\widehat{\widetilde{m}}(s,y)-\widehat{m}(0,x_{0}) =-i\alpha y\widehat{%
\widetilde{m}}(s,y)+\frac{1}{2}\beta ^{2}(iy)^{2}\widehat{\widetilde{m}}(s,y)
+\beta iy(\widehat{m}(.,y)\diamond \mathring{B}(.))(s),
\end{eqnarray*}
or
\[
\widetilde{\widehat{m}}(s,y)[s+i\alpha y+\frac{1}{2}\beta ^{2}y^{2}]=%
\widehat{m}(0,x_{0})+\beta iy(\widetilde{\widehat{m}(.,y)\diamond 
\mathring{B}(.)})(s),
\]
or
\begin{align}
\widetilde{\widehat{m}}(s,y)&=\frac{\widehat{m}(0,x_{0})}{s+i\alpha y+\frac{1
}{2}\beta ^{2}y^{2}}+\frac{\beta iy}{s+i\alpha y+\frac{1}{2}\beta ^{2}y^{2}}(
\widetilde{\widehat{m}(.,y)\diamond \mathring{B}(.)})(s)\nonumber\\
&=\frac{\widehat{m}(0,x_{0})}{s+i\alpha y+\frac{1
}{2}\beta ^{2}y^{2}}+\beta iy L(e^{(-i\alpha y-\frac{1}{2}\beta^{2}y^{2})t}(s)
\widetilde{\widehat{m}(.,y)\diamond \mathring{B}(.)})(s).
\end{align}
Put $g(t)=e^{(-i\alpha y-\frac{1}{2}\beta^{2}y^{2})t}$ and $h(t)= \widehat{m}(t,y)\diamond \mathring{B}(t)).$\\
Taking inverse Laplace transform, we get%
\begin{eqnarray}
\widehat{m}(t,y) &=&\widehat{m}(0,x_{0})\exp ((i\alpha y-\frac{1}{2}\beta
^{2}y^{2})t)  \label{(10)}+\beta iyL^{-1}(Lg \cdot Lh)(t,y)  \nonumber \\
&=&\widehat{m}(0,x_{0})\exp ((i\alpha y-\frac{1}{2}\beta ^{2}y^{2})t) 
\nonumber +\beta iy(g\ast h)(t,y),  \nonumber
\end{eqnarray}%
where
\[
(g\ast h)(t)=\int_{0}^{t}g(s)h(t-s)ds.
\]%
Recall that%
\begin{equation}
\int_{\mathbb{R}}e^{-ay^{2}-2by}dy=\sqrt{\frac{\pi }{a}}e^{\frac{b^{2}}{a}};\quad a>0.
\label{(11)}
\end{equation}%
Hence%
\begin{eqnarray*}
F^{-1}(g) &=&\frac{1}{2\pi }\int_{\mathbb{R}}e^{i\alpha yt-\frac{1}{2}\beta
^{2}y^{2}t+iyx}dy=\frac{1}{2\pi }\int_{\mathbb{R}}e^{-\frac{1}{2}\beta
^{2}ty^{2}-2yi(\frac{1}{2}\alpha t-\frac{1}{2}x)}dy \\
&=&\frac{1}{2\pi }\sqrt{\frac{\pi }{\frac{1}{2}\beta ^{2}t}}\exp (\frac{%
i^{2}(\frac{1}{2}\alpha t+\frac{1}{2}x)^{2}}{\frac{1}{2}\beta ^{2}t}) \\
&=&\frac{1}{\sqrt{2\pi \beta ^{2}t}}\exp (-\frac{(\alpha t+x)^{2}}{2\beta
^{2}t}) \\
&=&\frac{1}{\sqrt{2\pi \beta ^{2}t}}\exp (-\frac{x^{2}}{2\beta ^{2}t}-\frac{%
\alpha x}{\beta ^{2}})\exp (-\frac{\alpha ^{2}t}{2\beta ^{2}})=:k(t,x).
\end{eqnarray*}%
Therefore $g(t,x)=F[k(t,.)](y)]$ and (\ref{(10)}) can be written%
\begin{eqnarray*}
\widehat{m}(t,y) &=&\widehat{m}(0,x_{0})F[k](t,y) \\
&&+\beta iy\int_{0}^{t}F[k(t-s,.)](y)F[m(s,y)\diamond  \mathring{B}(s))]ds
\end{eqnarray*}%
Taking inverse Fourier transform we get, with $k^{\prime }=\frac{d}{dx}k(t,x)
$%
\begin{eqnarray*}
m(t,x) &=&m(0,x_{0})\ast k)(t,x) \\
&+&F^{-1}[\beta \int_{0}^{t}F[k^{\prime }(t-s,\cdot)](y)F[m(s,y)\diamond 
\mathring{B}(s))]](t,x) \\
&=&\int_{\mathbb{R}}\delta _{x_{0}}(x-y)k(t,y)dy \\
&&+\beta \int_{0}^{t}(\int_{\mathbb{R}}k^{\prime
}(t-s,x-y)m(s,y)\diamond \mathring{B}(s)dsdy \\
&=&k(t,x-x_{0})+\beta \int_{\mathbb{R}}(\int_{0}^{t}k^{\prime
}(t-s,x-y)m(s,y)dB(s))dy.
\end{eqnarray*}
We have proved the following:
\begin{theorem}\label{th5.5}
Suppose $\alpha$ and $\beta$ are constants and that the Donsker delta measure is absolutely continuous with respect to Lebesgue measure. Then the Donsker delta function $m(t,x)=\delta_{X(t)}(x)$ of the corresponding McKean-Vlasov process is a solution in $(\mathcal{S})^{*}$ of the following stochastic Volterra equation:
\begin{align*} \label{m}
m(t,x) &=k(t,x-x_{0})+\beta \int_{\mathbb{R}}(\int_{0}^{t}k^{\prime
}(t-s,x-y)m(s,y)dB(s))dy, 
\end{align*}
where
\begin{align*}
    k(t,z)=\frac{1}{\sqrt{2\pi \beta ^{2}t}}\exp (-\frac{z^{2}}{2\beta ^{2}t}-\frac{%
\alpha z}{\beta ^{2}})\exp (-\frac{\alpha ^{2}t}{2\beta ^{2}});\quad k^{\prime }(u,z)=\frac{d}{dz}k(u,z).
\end{align*}
\end{theorem}

\begin{remark}

If $\alpha =0,\beta =1$ we get%
\begin{eqnarray*}
k(t,z) &=&\frac{1}{\sqrt{2\pi t}}\exp (-\frac{z^{2}}{2t}) \\
k^{\prime }(u,z) &=&-\frac{z}{u}\frac{1}{\sqrt{2\pi u}}\exp (-\frac{z^{2}}{2u}).
\end{eqnarray*}
For comparison, recall that the density of Brownian motion at $t,x$ (when starting at $x_{0}$) is%
\[
p(t,x)=\frac{1}{\sqrt{2\pi t}}\exp (-\frac{(x-x_{0})^{2}}{2t}).
\]
\end{remark}

\section{Acknowledgments}
We are grateful to Frank Proske for helpful comments.

\end{document}